\documentclass[a4paper]{article}
\usepackage[english]{babel}
\usepackage{amsthm}  
\usepackage{amsmath} 
\usepackage{amssymb} 
\usepackage{mathabx} 
\usepackage{graphicx}
\usepackage{fancyhdr}
\usepackage{setspace}
\usepackage{lineno}	 
\usepackage{lastpage}
\usepackage{enumitem}
\usepackage{mathtools}
\usepackage{epsfig} 
\usepackage[outdir=./]{epstopdf} 
\usepackage[margin=1in]{geometry}
\usepackage{caption}  
\usepackage[round]{natbib}
\usepackage{tocloft}	
\allowdisplaybreaks

\newcommand{\mytitle}{Adaptive Learning of Hybrid Models for Nonlinear Model Predictive Control of Distillation Columns}
\newcommand{\myshorttitle}{Adaptive Learning}
\newcommand{\myauthor}{Jannik T. L\"uthje, Jan C. Schulze, Adrian Caspari, Adel Mhamdi, Alexander Mitsos, and Pascal Sch\"afer$^{*}$} 
\newcommand{\myauthorshort}{L\"uthje et al.}
\author{\myauthor}
\usepackage[colorlinks,linkcolor=blue,citecolor=blue,urlcolor=blue,pdftitle={\mytitle},pdfauthor={Pascal Sch\"afer and Alexander Mitsos}]{hyperref} 

{
	\fancyhead[L]{\myshorttitle}
	\fancyhead[R]{\small{\textit{\the\day.\the\month.\the\year}}}
	\fancyfoot[L]{\copyright \small{\textit{\myauthorshort}}}
	\fancyfoot[R]{Page \thepage\ of \pageref*{LastPage}}
	\cfoot{\textit{}}
}

\fancypagestyle{firststyle}
{
	 \fancyhead[L]{\copyright \small{\textit{\myauthorshort}}}
	 \fancyhead[C]{\textit{}}
   \fancyhead[R]{Page \thepage\ of \pageref*{LastPage}}
	 \fancyfoot[L]{\footnotesize \rule{0.25\textwidth}{0.4pt} \newline Corresponding author: $^*$P. Sch\"afer \newline AVT Process Systems Engineering, RWTH Aachen University, 52056 Aachen, Germany \newline E-mail: pascal.schaefer@avt.rwth-aachen.de}
   \fancyfoot[C]{}
	 \fancyfoot[R]{}
}

\theoremstyle{remark} 

\makeatletter
\let\@addpunct\@gobble
\g@addto@macro{\thm@space@setup}{\thm@headpunct{}} 
\renewcommand\@biblabel[1]{#1.}
\makeatother

\renewenvironment{abstract}{\noindent\textbf{Abstract:}}{}

\newenvironment{acknowledgements}{\noindent\footnotesize\textbf{Acknowledgements}}{}

\begin{document}

\thispagestyle{firststyle}
\begin{flushleft}\begin{large}\textbf{\mytitle}\end{large} \end{flushleft}
\myauthor 

\begin{flushleft}\begin{small}
RWTH Aachen University\\
AVT - Aachener Verfahrenstechnik\\
Process Systems Engineering\\
52072 Aachen, Germany\\[0.25cm]

\end{small}
\end{flushleft}

\noindent\hrulefill
\vspace{0.25cm} \\
\begin{abstract}
Nonlinear model predictive control (NMPC) requires accurate and computationally efficient plant models. Our previous work has shown that the classical compartmentalization model reduction approach for distillation columns can be enhanced by replacing parts of the system of equations by artificial neural networks (ANNs) trained on offline solved solutions to improve computational performance. In real-life applications, the absence of a high-fidelity model for data generation can, however, prevent the deployment of this approach. Therefore, we propose a method that utilizes solely plant measurement data, starting from a small initial data set and then continuously adapting to newly measured data. The efficacy of the proposed approach is examined in silico for a distillation column from literature. To this end, we first adjust our reduced hybrid mechanistic/data-driven modeling approach that originally builds on compartmentalization to a stage-aggregation procedure, tailoring it for the application within the adaptive learning framework. Second, we apply an adaptive learning algorithm that trains the ANNs replacing the stationary stage-to-stage calculations on newly available data. We apply the adaptive learning of the hybrid model within a regulatory NMPC framework and conduct closed-loop simulations. We demonstrate that by using the proposed method, the control performance can be steadily improved over time compared to a non-adaptive approach while being real-time applicable. Moreover, we show that the performance when using either a model trained on excessive amounts of offline generated data or the original high-fidelity model can be approached in the limit.\\
\end{abstract}

\vspace{0.25cm}
\noindent\hrulefill

\section{Introduction}
Nonlinear model predictive control (NMPC) requires accurate dynamic process models resulting in optimization problems that are typically not solvable in real time. This is particularly the case for distillation columns if applying rigorous full-order stagewise models. Thus, several reduction methods for dynamic models have been developed. A short review of these methods for distillation columns is given in \cite{Schafer.2020}. These methods are the collocation approach (\cite{Cho.1983}, applied in \cite{Cao.2016}), the compartmentalization approach (\cite{Benallou.1986}, applied in \cite{Schafer.2019}), and the wave propagation approach (\cite{Marquardt.1989}, applied in \cite{Caspari.2020}). For a broader review of model reduction methods we refer to \cite{WMarquardt2002}.\\

\noindent The computational benefit introduced by the model reduction can be further extended by introducing machine learning elements (cf. \cite{Schafer.2020}). In particular, we have shown in previous work that compartmentalized models can be substantially enhanced by substituting stationary stage-to-stage calculations by artificial neural networks (\cite{Schafer.2019}). In a subsequent work, we demonstrate that this method enables real-time capable NMPC (\cite{Schafer2019b}). Until now, these hybrid mechanistic/data-driven modeling approaches, however, rely on the availability of large data sets before operation. These data sets are mostly generated offline, utilizing the full-order stagewise model. However, the existence of a full-order model is not guaranteed in real-life application. Thus, recent work has focused on the on-line adaption of the process model to account for an inaccurate initial model for operational purposes (e.g., \cite{Tsay.2020}).\\

\noindent In this work, we propose a procedure to continuously adapt our published hybrid model (\cite{Schafer.2019}) to new measured data, thus, removing the requirement of the existence of a full-order model. To this end, we reformulate the hybrid compartmentalized model as a hybrid stage-aggregation model (cf. \cite{Linhart.2010}), which allows for a stricter separation into steady-state blocks to be replaced by artificial neural networks (ANNs) and dynamic single stages with measurable states. As the data-driven model parts have to be continuously adapted to new data, we utilize an on-line learning algorithm that is repeatedly applied, as opposed to a batch learning algorithm that is only applied before operation. We utilize this idea in and in silico control case study from literature and compare to a hybrid control scheme based on offline-generated data and an ideal NMPC using the exact full-order stagewise model.\\

\noindent The remainder of this work is structured as follows: first, we describe the proposed hybrid stage-aggregation model and the implemented adaptive learning algorithm. Next, the case study is presented. Afterwards, we investigate the control performance and real-time applicability of the proposed approach. Finally, an outlook on further work is given.

\section{Method}
\subsection{Hybrid stage-aggregation model}
The full-order model, which we consider as basis for this work, is taken from \cite{Rehm.1996}, where the model is formulated for a distillation column that is described in \cite{Allgower.1992}. It assumes constant molar holdup, ideal thermodynamics, no pressure losses, and no hydrodynamic resistance. These assumptions reduce the model to one mass balance and one thermodynamic equation for each stage (cf. Appendix \ref{sec:full}).\\

\noindent In a compartmentalization approach (cf. \cite{Benallou.1986}), the distillation column is split into multiple compartments. We then assume that the single-stage dynamic behavior can be neglected compared to the overall dynamic behavior of the entire compartment.
The final system of equations thus consists of the dynamic compartment balances and steady-state equations for each stage inside the compartment. One stage per compartment does not have a steady-state equation and is referred to as \textit{sensitivity stage}. Thereby, this type of model reduction does not change the steady-state behavior of the column; however, it affects the dynamic response as well as the stiffness of the differential-algebraic system.\\

\noindent For our adaptive approach, we use the stage-aggregation analogy presented by \cite{Linhart.2010}, reformulating the compartmentalized model in an exact way. That is, the stages of the column are split into two groups, aggregation stages and non-aggregation stages.
Aggregation stages have increased holdup, whereas the holdup of non-aggregation stages is set to zero, making them steady-state.
The advantage of this approach is that aggregation stages correspond to actual stages of the column, allowing to identify the states of the reduced model from real-world measurements. Formally, this procedure corresponds to multiplying the holdups of aggregation stages with a holdup factor \(H\) and reducing the holdup of steady-state trays to zero (cf. Appendix \ref{sec:reduced}). Feed stages, the condenser, and the reboiler have to be aggregation stages.\\
\newline
\begin{figure}[tb]
 \centering
 \includegraphics[width=0.75\linewidth]{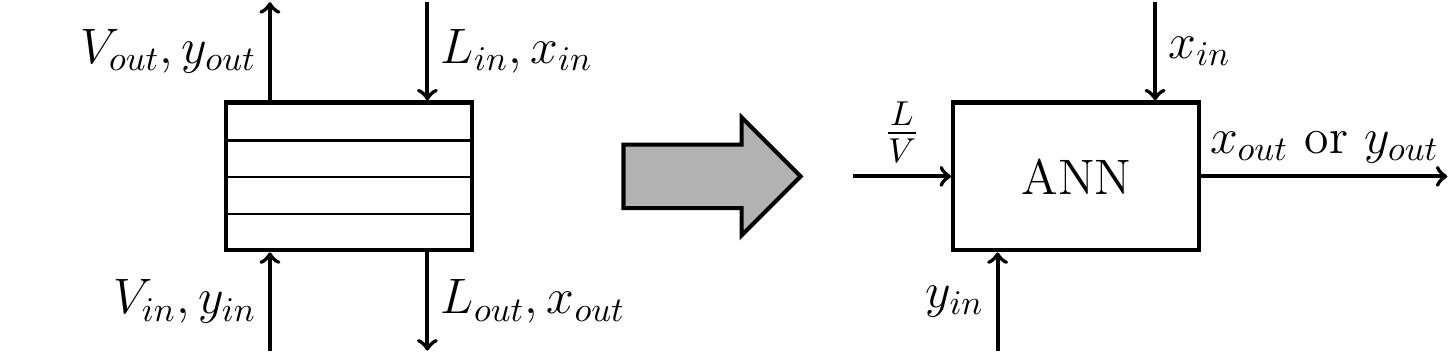}
 \caption{Replacement of a steady-state column section with an ANN.}\label{fig:Hybridisation}
\end{figure}
\noindent Next, we reduce the total number of equations of the reduced model by replacing sections of steady-state trays, i.e., stationary column sections between two consecutive aggregation stages, with ANNs (Figure \ref{fig:Hybridisation}), adapting the approach from our previous work considering compartmentalization (\cite{Schafer.2019}) to stage-aggregation. Note that in contrast to hybrid compartment models, the differential states of hybrid stage-aggregation models are not direct inputs of the surrogate model replacing the stationary parts, leading to generally lower input dimensionalities. Each column section maps the out-going liquid stream of the upper aggregation stage and the out-going vapor stream of the lower aggregation stage to the in-going vapor stream of the upper aggregation stage and the in-going liquid stream of the lower aggregation stage. Due to the simplicity of the assumed full-order model, the liquid and vapor flow, \(L\) and \(V\), are constant over the column section. Thus, the ANN does not need to consider the total molar flows as outputs. From the remaining two output concentrations (\(x_{out}\) and \(y_{out}\)), one can be computed by a mass balance around the entire steady-state column section. Since all trays to be replaced have no holdup, only the ratio of flows, \(\frac{L}{V}\), needs to be considered. Thus, the ANNs replacing the column sections have three inputs and one output. Logarithmic scaling is applied to the molar fractions as proposed by \cite{Cao.2016}.
\begin{figure*}[tb]
 \centering
 \includegraphics[width=\textwidth]{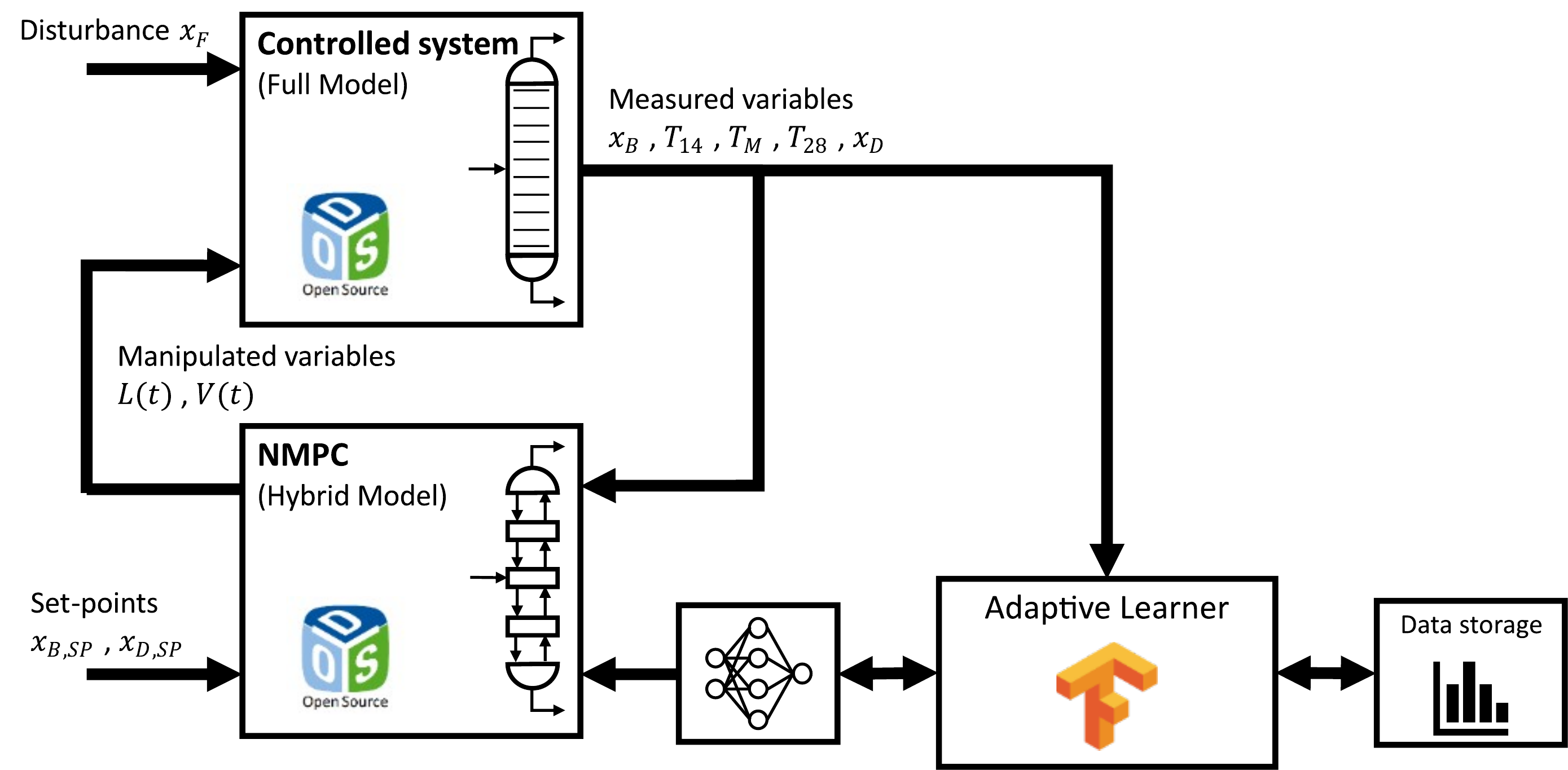}
 \caption{Implementation of the adaptive control scheme.}\label{fig:Implementation}
\end{figure*}

\subsection{Adaptive learning algorithm}
There are different interpretations of the definition of on-line learning in literature (cf. \cite{PerezSanchez.2018}), in this work we focus on learning of a continuous data stream.
As opposed to batch learning, the complexity of the data is not known beforehand. To ensure an adequate model size at every time step, constructive algorithms are proposed in literature (e.g., \cite{Ma.2003}). Constructive algorithms start with a small model size and continuously increase the model size whenever the current model size shows to be incapable of representing the new data.\\

\noindent We implement an adaptive learning algorithm based on the work of \cite{Chen.1993}.
The overall objective of the algorithm is to improve performance on newly presented data while maintaining performance on previously seen data.
For this purpose, all seen data is stored. In each training cycle, a sample of the old data is chosen and mixed with the entire new data set.
To select data from the storage, latin hypercube sampling (LHS) is performed in the data region of the data storage. For each sampled point, the nearest neighbor in the data storage is determined. These neighbors finally become part of the training data.
With this approach, we build a training data set, which is as equally distributed as possible, out of the unequally distributed data storage.
Training of the ANNs is performed with the Levenberg-Marquardt algorithm.
The current model is used as a starting point.
The balance between improving performance on new data and maintaining performance on previously seen data is tuned by adding a factor to the weights of new data points.\\

\noindent As described earlier, a constructive model building approach is applied to ensure adequate model size at all times.
That is, whenever a performance goal cannot be met after the first training cycle, the model size is increased.
To this end, the model is restored to its state before the failed training cycle. Then, a new node is added to the single hidden layer.
To initialize this new node, a training cycle is performed in which only the weights corresponding to the new node are changed.
Considering the linear activation function of the output layer, this is equivalent to fitting a one-node ANN to the prediction error of the old model on the training set. After the new node is initialized, another training cycle is performed using the enlarged model as an initial point, this time adjusting all weights. If the performance goal is still not reached, the procedure is repeated.

\subsection{Implementation and case study}
The described framework is applied to a distillation column commonly used in literature, first and in most detail described in \cite{Allgower.1992}. The column consists of 40 stages and separates a mixture of methanol and propanol. The control case study is taken from \cite{Diehl.2002}, meaning that the feed concentration is the unmeasured disturbance. We apply set-points on the molar fractions of both column products.
The control objective \(\varphi\) is the integrated sum of the set-point deviations:
\begin{equation}
 \varphi = \int_0^{T_\text{P}} \left(x_{\text{B,SP}} - x_\text{B}\left(t\right)\right)^2 + \left(x_{\text{D,SP}} - x_\text{D}\left(t\right)\right)^2\; \text{d} t\; .
\end{equation}
\begin{table}[tb]
 \centering
 \caption{MPC settings}\label{tab:MPC_Settings}
 \begin{tabular}{l c r}
  \hline
  Setting            & Symbol &  Value \\
  \hline
  Control horizon    & \(T_\text{C}\) & 600s\\
  Prediction horizon & \(T_\text{P}\) & 1200s \\
  Control intervals  & \(N\) & 10    \\
  Sampling time      & \(T_\text{s}\) & 60s    \\
  \hline
 \end{tabular}
\end{table}

\noindent The MPC settings, summarized in Table \ref{tab:MPC_Settings}, are also taken from \cite{Diehl.2002}, except for the sampling time which is increased from 10s to 60s to allow for a longer solution time of the dynamic optimization problem but still appears appropriate considering the inertia of the column.\\

\noindent Concerning the measurements, \cite{Allgower.1992} specify temperature measurements on stages 14 and 28. As the measurements allow to determine the state of these stages, they are selected as aggregation stages. Additionally, we assume a temperature measurement for the feed stage. Finally, we assume concentration measurements for the top and bottom product for the sake of simplicity. Thus, the state of all aggregation stages can be directly derived from the measurements. In a compartmentalization approach this direct identification of the differential variables would not work as compartments do not correspond to actual column trays. In cases where temperatures are measured only, suitable state estimation techniques can be applied to determine the states of the reduced model.\\

\noindent Training data for the four ANNs can be calculated using mass balances around the aggregation stages. An estimate of the feed flow molar fraction \(x_F\) can be determined by a mass balance around the entire hybrid model of the column. Note, that both of these procedures rely on the underlying assumption that the hybrid model matches the full-order model exactly. However, assuming sufficient accuracy of the ANNs, this is only true in the steady-state case. For non-steady-state operation, an error will be made due to lower-order approximation of the dynamics. To avoid learning wrong data, a weight is added to each data point, which diminishes as the column moves away from steady-state operation.\\

\noindent A closed-loop control framework (cf. Figure \ref{fig:Implementation}) is set up to assess the proposed method. For this, the controlled system is simulated using the full-order stagewise model as plant replacement. The NMPC employs the proposed hybrid stage-aggregation model as controller model. Both the plant replacement and the controller model are implemented in Modelica and made accessible for our in-house software-package DyOS (\cite{Caspari.2019}) via export as Functional Mockup Unit (FMU). DyOS enables access to LIMEX (\cite{Schlegel.2004}) for state and sensitivity integration and SNOPT (\cite{Gill.2005}) for solving the dynamic optimization problems in a single-shooting approach.\\

\noindent The adaptive learner reads the measurement data generated by the simulated column, derives the training data, and adapts the ANNs of the hybrid model. The algorithm is implemented in Python using the Tensorflow library (cf. \cite{Abadi.2015}). The interconnections of the described components are also managed in Python. All calculations are performed on an Intel\copyright\ Xeon\copyright\ Gold 5117 CPU. The Intel\copyright\ Optimization for Tensorflow is used.\\

\noindent The adaptive learning algorithm and the NMPC require an initial model to function. For this, we generate open-loop data by performing step tests on the plant replacement. For the step-tests, we combine individual step tests on the manipulated variables \(L\) and \(V\), and the disturbance \(x_F\). The height of these steps and their time are determined using latin hypercube sampling. Overall, 20 open-loop trajectories are performed to build the initial data set. These trajectories are then transformed to training data onto which the ANNs are trained.\\

\noindent Overall, we compare four different approaches in this work:
\begin{enumerate}[label=(\roman*)]
  \item The adaptive approach proposed in this work.
  \item A control scheme based only on the initial model used in the adaptive approach.
  \item An approach using an offline-trained model from sampled data using the full-order stagewise model as done in previous work. Here, a larger, broader, and more equally distributed data set is used. Also, the data has higher quality, as it is not influenced by the dynamic mismatch between the full-order and the reduced model.
  \item As a benchmark, we apply the ideal NMPC with full-state feedback and without plant-model mismatch, i.e., the controller model is exact.
\end{enumerate}

\section{Results}
\subsection{Control performance}

\begin{figure}[tb]
  \centering
 \includegraphics[width=0.6\linewidth]{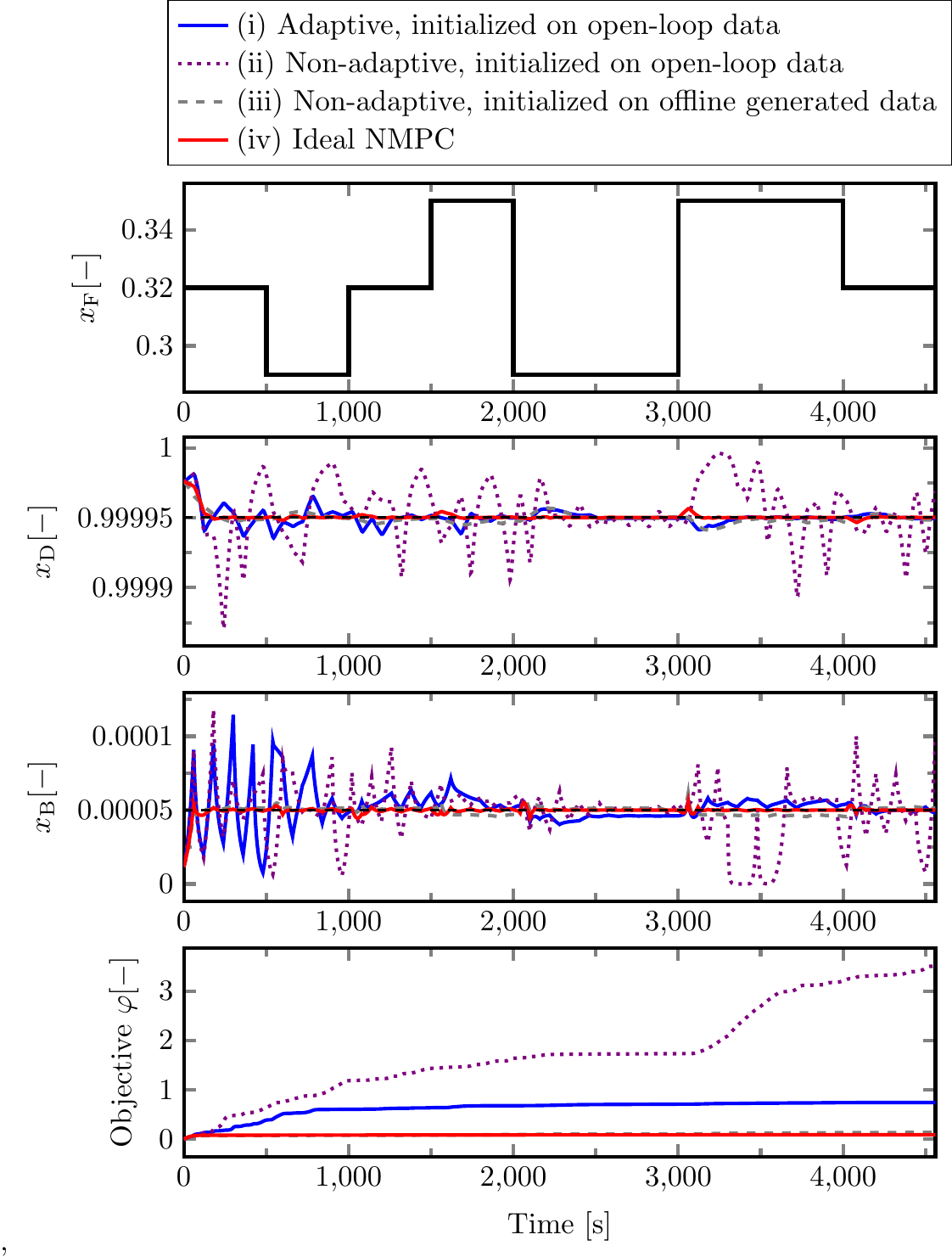}
 \caption{Comparison of four different control schemes. The upper plot shows the trajectory of the unmeasured disturbance \(x_\text{F}\). The remaining plots show trajectories of the plant replacement simulated in closed-loop for all approaches. These are the condenser molar fraction \(x_\text{D}\), the reboiler molar fraction \(x_\text{B}\), and the NMPC objective \(\varphi\).  The set-points on \(x_\text{D}\) and \(x_\text{B}\) are 0.99995 and 0.00005 respectively.}\label{fig:Results}
\end{figure}

\noindent Figure \ref{fig:Results} presents the results of the four previously described approaches. The non-adaptive approach based only on open-loop data shows strongly oscillating behavior throughout the entire horizon. Large deviations from the set-points can be observed (e.g., the reboiler concentration around \(t=3500\text{s}\)). The introduction of the adaptive learning algorithm improves the control performance substantially. After \(t=1000\text{s}\), almost no oscillating behavior can be observed. Looking at the objective, the integrated sum of squared errors of the set-point deviations, we see a rapid flattening of the trajectory. At the end of the case study, the adaptive approach reaches a control performance comparable to the approach utilizing excessive amounts of offline-generated data. Finally, we remark that the ideal NMPC without plant/model-mismatch allows for only slight improvements over the approach using offline-trained models and thus also over the final performance of the proposed adaptive control scheme.

\subsection{Real-time applicability}

\begin{figure}[tb]
  \centering
 \includegraphics[width=0.6\linewidth]{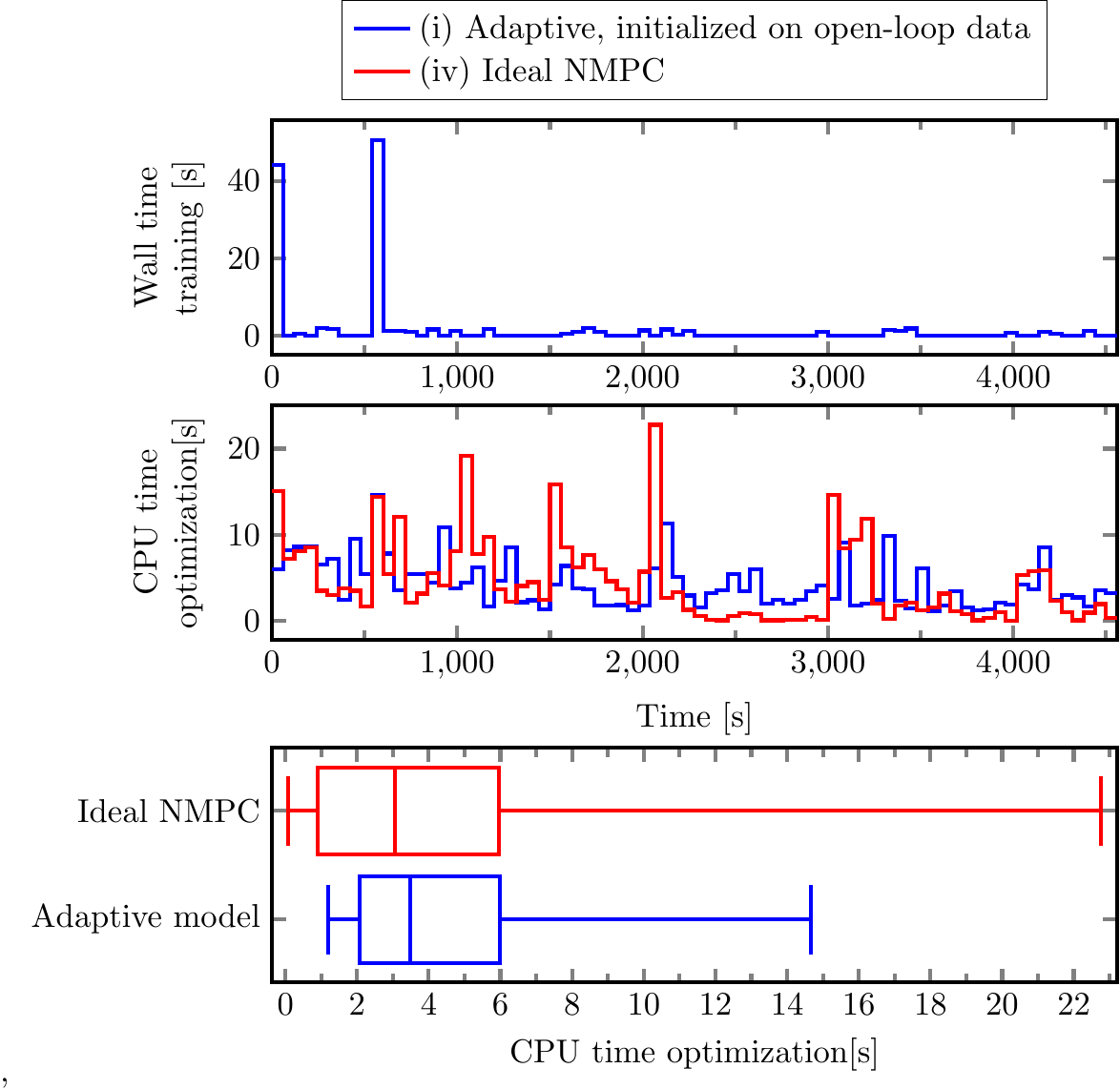}
 \caption{Overview of computational requirements for the adaptive learning approach and the ideal NMPC. The top plot shows the required time to adapt all ANNs of the adaptive hybrid model. Next, the times required to solve the dynamic optimization problems are presented. The last plot is a box plot comparing the solution times of the dynamic optimizations for the adaptive approach and the ideal NMPC. The box plot visualizes the minimum, the 25\textsuperscript{th} percentile, the median, the 75\textsuperscript{th} percentile, and the maximum of the data.}\label{fig:CPU}
\end{figure}

\noindent Besides the control performance, the computational time requirement is another important aspect when rating control schemes, as real-time applicability has to be maintained. For the adaptive approach, the time required to solve the dynamic optimization problem and the time required for adapting the ANNs have to be considered. Note, that the training of the ANNs can be performed in parallel which is not implemented at the moment.\\

\noindent Figure \ref{fig:CPU} gives the computational times. The training time shows two large peaks with \(t_{Train} > 40\text{s}\) at \(t = 0\text{s}\) and \(t = 500\text{s}\). At these times large ANN size increases take place. Besides these peaks, the training time stays below \(5\text{s}\). The dynamic optimization of the adaptive approach finishes in under \(15\text{s}\) throughout the entire case study. The total computational time (sum of training and optimization) only slightly exceeds the sampling time of \(T_\text{s} = 60\text{s}\) once (\(t = 500\text{s}\)), demonstrating the real-time applicability of the approach in this setting. Concerning the computational time for the dynamic optimization, we also compare the adaptive approach with the ideal NMPC. To this end, one has to distinguish between two different control scenarios with the first being right after a step in which the disturbance has occurred (cf. \(x_\text{F}\) plot in Figure \ref{fig:Results}). Here, a high computational effort is required to solve the dynamic optimization problem as the solution of the previous time step is an inadequate initial guess. For this scenario, we see a large computational benefit of the hybrid model compared to the full-order model used in the ideal NMPC (e.g., at \(2000\text{s}\)), as both approaches require a large number of iterations to solve, thus, allowing the hybrid approach to fully show its computational benefits. In contrast, the second scenario considers times when there is no change in the disturbance.
In that case, the lack of plant/model-mismatch of the ideal NMPC and, thus, the existence of better initial guesses leads to a computational advantage for the ideal NMPC during these times, as it requires less iterations to improve the solution. This results in the ideal NMPC having a broader distribution of solution times than the adaptive approach (cf. box plot in Figure \ref{fig:CPU}). Note however that the high computational times after occurrence of disturbances are the more crucial ones with regard to real-time applicability. Furthermore and following our previous results (cf. \cite{Schafer2019b}), we suspect computational time savings when applying the hybrid model in the adaptive approach to become even more significant if the control task becomes more complex (e.g., more complex model, more severe disturbances, or economic objectives).

\section{Conclusion}

We develop an adaptive control scheme for distillation columns that continuously trains the data-driven parts of a reduced hybrid controller model based on measurement data. For this purpose, we present a hybridized stage-aggregation approach and implement an adaptive learning algorithm. The approach manages to improve control performance over time, compared to a non-adaptive approach, while maintaining real-time applicability.\\

\noindent In this work, the method is demonstrated in a single-column case study based on a simplified full-order model. Introducing more complexity to the model would result in larger ANNs (e.g., non-constant flows would require at least one additional ANN output).
As for the framework, the robustness of the approach could be improved by adding constraints that restrict the dynamic optimization to previously explored data regions. Furthermore, different approaches to deal with the unequally distributed closed-loop data can be investigated.
Also, the approach can be applied to the use case of a drift in the underlying process. By implementing a forgetting factor for stored data points, the data set could be replaced over time to match the new process, thus adapting to the new and different process model. \\

\begin{acknowledgements}
The authors gratefully acknowledge the financial support of the Koper\-nikus project Syn\-Ergie by the Federal Ministry of Education and Research (BMBF) and the project supervision by the project management organization Projekttr\"ager J\"ulich.
\end{acknowledgements}

\bibliographystyle{abbrvnat}       
\bibliography{ifacconf}   

\appendix
\numberwithin{equation}{section}
\section{Full-order model}\label{sec:full}
Standard tray:
\begin{equation}
 n_i \frac{dx_i}{dt} = L^*\left(x_{i,in}-x_i\right) + V\left(y_{i,in}-y_i\right)
\end{equation}
with \(L^*=L\) in the rectifying section and \(L^*=L+F\) in the stripping section.\\
\newline
Condenser (\(D\)):
\begin{equation}
 n_D \frac{dx_D}{dt} = V\left(y_{in}-x_D\right)
\end{equation}
Feed tray (\(M\)):
\begin{equation}
n_M \frac{dx_M}{dt} = L\left(x_{in} - x_M\right) + V\left(y_{in}-y_M\right) + F\left(x_F-x_M\right)
\end{equation}
Reboiler (\(B\)):
\begin{equation}
 n_B \frac{dx_B}{dt} = \left(L+F\right)\left(x_{in}-x_B\right) + V\left(x_B-y_B\right)
\end{equation}
Thermodynamics:
\begin{equation}
 y_i = \frac{\alpha\ x_i}{1+\left(\alpha - 1\right) x_i}
\end{equation}

\section{Reduced model}\label{sec:reduced}
Aggregation stage:
\begin{equation}
 H_i n_i \frac{dx_i}{dt} = L^*\left(x_{i,in}-x_i\right) + V\left(y_{i,in}-y_i\right)
\end{equation}
Non-aggregation stage:
\begin{equation}
 0 = L^*\left(x_{i,in}-x_i\right) + V\left(y_{i,in}-y_i\right)
\end{equation}
Condenser (\(D\)):
\begin{equation}
 H_D n_D \frac{dx_D}{dt} = V\left(y_{in}-x_D\right)
\end{equation}
Feed tray (\(M\)):
\begin{equation}
H_M n_M \frac{dx_M}{dt} = L\left(x_{in} - x_M\right) + V\left(y_{in}-y_M\right) + F\left(x_F-x_M\right)
\end{equation}
Reboiler (\(B\)):
\begin{equation}
 H_B n_B \frac{dx_B}{dt} = \left(L+F\right)\left(x_{in}-x_B\right) + V\left(x_B-y_B\right)
\end{equation}
\end{document}